\documentclass[12pt]{amsart}

\usepackage{enumerate, amsmath, amsthm, amsfonts, amssymb, xy,  mathrsfs, graphicx, paralist, fancyvrb, eucal}
\usepackage[usenames, dvipsnames]{xcolor}
\usepackage[margin=1in]{geometry} 
\usepackage[bookmarks, colorlinks=true, linkcolor=blue, citecolor=blue, urlcolor=blue]{hyperref}

\input xy
\xyoption{all}

\numberwithin{equation}{section}
\newtheorem{theorem}[equation]{Theorem}

\newtheorem{proposition}[equation]{Proposition}

\theoremstyle{definition}
\newtheorem{rmk}[equation]{Remark}
\newenvironment{remark}[1][]{\begin{rmk}[#1] \pushQED{\qed}}{\popQED \end{rmk}}
\newtheorem{eg}[equation]{Example}

\newtheorem{defn}[equation]{Definition}
\newenvironment{definition}[1][]{\begin{defn}[#1]\pushQED{\qed}}{\popQED \end{defn}}

\newcommand{\bC}{\mathbf{C}}
\newcommand{\cC}{\mathcal{C}}

\newcommand{\rH}{\mathrm{H}}

\newcommand{\bN}{\mathbf{N}}

\newcommand{\bP}{\mathbf{P}}

\newcommand{\cV}{\mathcal{V}}

\newcommand{\bZ}{\mathbf{Z}}

\newcommand{\bk}{\mathbf{k}}

\renewcommand{\phi}{\varphi}
\renewcommand{\emptyset}{\varnothing}

\newcommand{\arxiv}[1]{\href{http://arxiv.org/abs/#1}{{\tt arXiv:#1}}}

\makeatletter
\def\Ddots{\mathinner{\mkern1mu\raise\p@
\vbox{\kern7\p@\hbox{.}}\mkern2mu
\raise4\p@\hbox{.}\mkern2mu\raise7\p@\hbox{.}\mkern1mu}}
\makeatother

\renewcommand{\hom}{\operatorname{Hom}}

\DeclareMathOperator{\Sym}{Sym}

\DeclareMathOperator{\Mod}{Mod}

\newcommand{\GL}{\mathbf{GL}}

\newcommand{\Conf}{{\rm Conf}}
\newcommand{\Sec}{{\rm Sec}}
\newcommand{\FI}{\mathbf{FI}}
\newcommand{\OI}{\mathbf{OI}}
\newcommand{\defi}[1]{{\bf \textsf{#1}}}
\newcommand{\multi}[2]{(\bZ^{#1}_{\ge 0})_{#2}}

\title{Noetherian properties in representation theory}

\author{Steven V Sam}
\address{Department of Mathematics, University of Wisconsin, Madison, WI}
\email{svs@math.wisc.edu}
\urladdr{\url{http://math.wisc.edu/~svs/}}
\thanks{SS was partially supported by NSF grant DMS-1500069.}
\date{July 3, 2017}

\begin{document}

\maketitle

\begin{abstract}
The goal of this expository article, based on a lecture I gave at the 2016 ICRA, is to explain some recent applications of ``categorical symmetries'' in topology and algebraic geometry with an eye toward twisted commutative algebras as a unifying framework. The general idea is to find an action of a category on the object of interest, prove some niceness property, like finite generation, and then deduce consequences from the general properties of the category. The key in these cases is to prove that the representation theory of this category is locally noetherian, and we will discuss an outline for such proofs.
\end{abstract}

\section{Introduction}

In the past few years, a number of stabilization phenomena, which uses representation theory in an essential way, have been observed in many different areas, including topology, algebraic geometry, group theory, and combinatorics. One striking feature is that these phenomena have been discovered independently by several groups of researchers working in these different areas, and while the language used by each is different, there are many underlying commonalities. Notably, one finds what we will call a ``categorical symmetry'' of the objects in question, and uses finite generation of this symmetry to deduce the desired stabilization of uniformity results.

The goal of this article is to present two such examples (one from topology and one from algebraic geometry) and to explain a possible unifying framework (twisted commutative algebras, aka tca's). While tca's are rather complicated objects at first glance, we hope that this presentation will motivate their study. We end with a brief discussion on future work that needs to be done.

This is by no means meant to be a survey of the subject. For that, one might consult \cite{draisma-notes} or \cite{farb-icm}. Instead, we aim to get the reader as quickly as possible to the philosophy behind ``representation stability'' and to get a feeling for how some of the results can be proven in some concrete examples.

\subsection*{Acknowledgements} 

This paper is based on a talk given at the 2016 International Conference on Representations of Algebras (ICRA) at Syracuse University, and I thank the organizers for the opportunity to present. I also thank Eric Ramos and Andrew Snowden for reading a draft of this article.

\section{Two examples}

\subsection{Cohomology of configuration spaces} \label{sec:conf-eg}

Let $X$ be a ``nice'' (connected, orientable, and with the homotopy type of a finite CW complex) manifold of dimension $\ge 2$. For $n \ge 0$, define the $n$th \defi{configuration space} by
\[
\Conf_n(X) = \{(x_1, \dots, x_n) \in X^n \mid x_i \ne x_j \text{ for } i \ne j\}.
\]
Let $\bk$ be a field. The following result can be found in \cite{church} (when ${\rm char}(\bk) = 0$) and \cite{CEFN} (in general).

\begin{theorem}[Church, Church--Ellenberg--Farb--Nagpal] \label{thm:CEFN}
With the notation above, the function
\[
n \mapsto \dim_\bk \rH^i(\Conf_n(X); \bk)
\]
is given by a polynomial function in $n$ for $n \gg 0$.
\end{theorem}

The symmetric group $S_n$ acts freely on $\Conf_n(X)$, and one can also study the cohomology of the quotient, though the behavior depends on the characteristic. The following result can be found in \cite{church} (when ${\rm char}(\bk) = 0$) and \cite{nagpal} (when ${\rm char}(\bk) > 0$).

\begin{theorem}[Church, Nagpal] \label{thm:conf-quotient} Keep the notation above.
  \begin{itemize}
  \item When ${\rm char}(\bk)=0$, $\dim_\bk \rH^i(\Conf_n(X)/S_n; \bk)$ is constant for $n \gg 0$.
  \item When ${\rm char}(\bk)=p>0$, $\dim_\bk \rH^i(\Conf_n(X)/S_n; \bk)$ is periodic with period a power of $p$ for $n \gg 0$.
  \end{itemize}
\end{theorem}

\begin{remark}
These results can be extended to more general classes of spaces, see \cite{petersen} and \cite{tosteson}.
\end{remark}

\subsection{Syzygies of secants of Veronese embeddings}

Let $\bk$ be a field of characteristic $0$. Let $B = \bigoplus_{n \ge 0} B_n$ be a $\bZ$-graded $\bk$-algebras with $B_0 = \bk$, which is generated in degree $1$, and satisfies $\dim_\bk B_1 < \infty$. Define the $d$th Veronese subring to be 
\[
B^{(d)} = \bigoplus_{n \ge 0} B_{nd}.
\]
This is also a $\bZ$-graded $\bk$-algebra which is generated in degree $1$ and $\dim_\bk  B^{(d)}_1 = \dim_\bk B_d < \infty$. In particular, we have a surjection $\Sym(B_d) \to B^{(d)}$ where $\Sym$ denotes the symmetric algebra. Denote the kernel by $I_B^{(d)}$. The general philosophy is that $B^{(d)}$ becomes better behaved with respect to various homological invariants as $d \to \infty$ (see \cite{backelin} and \cite{ERT} for examples related to Koszulity).

The tensor product $\Sym(V) \otimes \Sym(V)$ is also an algebra via the multiplication 
\[
(p_1 \otimes p_2)(q_1 \otimes q_2) = p_1q_1 \otimes p_2q_2
\]
on simple tensors (and extended linearly). We define a map 
\[
\Delta \colon \Sym(V) \to \Sym(V) \otimes \Sym(V)
\]
on the degree $1$ piece $V$ by $v \mapsto v \otimes 1 + 1 \otimes v$ and extend it to higher degrees by requiring that $\Delta$ is a ring homomorphism. Given ideals $I, J \subset \Sym(V)$, their \defi{join} is the ideal $I \star J$, which is the kernel of the composition 
\[
\Sym(V) \xrightarrow{\Delta} \Sym(V) \otimes \Sym(V) \to \Sym(V)/I \otimes \Sym(V)/J.
\]
Finally, define the \defi{secant ideals} of $I$ by $\Sec^1 I = I$ and $\Sec^r I = \Sec^{r-1} I \star I$. The following is the main result of \cite{secver}.

\begin{theorem}[Sam] \label{thm:secver}
Fix $B$ as above, and an integer $r \ge 1$. There exists a constant $C_B(r)$ such that for all $d$, the ideal $\Sec^r I_B^{(d)}$ is generated in degrees $\le C_B(r)$.
\end{theorem}

For a geometric interpretation of this statement, consider the cases when $I,J \subseteq \Sym(V)$ are the homogeneous prime ideals defining projective varieties $V(I), V(J)$ in projective space $\bP(V^*)$. Given a point $v$ in the vector space $V^*$, let $[v] \in \bP(V^*)$ denote the line it spans. Then
\[
V(I\star J) = \overline{ \{ [x + y] \mid [x] \in V(I), [y] \in V(J) \}}.
\]
To be more precise, we are taking the set of all lines spanned by elements of the form $x+y$ where the line spanned by $x$ is in $V(I)$ and the line spanned by $y$ is in $V(J)$, and the overline denotes Zariski closure. Also, if $B$ is the homogeneous coordinate ring of $X \subset \bP(V^*)$, then $V(I_B^{(d)}) \subset \bP((\Sym^d V)^*)$ is the image of $X$ under the $d$th Veronese map (the $d$th power map).

Phrased this way, Theorem~\ref{thm:secver} says that given an embedded projective variety $X \subseteq \bP^n$, there is a bound $C_X(r)$ so that the full ideal of the $r$th secant variety of any Veronese re-embedding of $X$ is generated in degree $\le C_X(r)$. 

\begin{remark}
The above leads to a relaxation of the problem: instead of finding degree bounds for the generators of the full ideal of the $r$th secant variety, look for degree bounds for generators of any ideal whose radical is the full ideal of the $r$th secant variety (i.e., just find equations that define the same set as the $r$th secant variety). In fact, analogous theorems for this set-theoretic relaxation exist for the $r$th secant variety of arbitrary Segre products of projective spaces \cite{draisma-kuttler} and for the $r$th secant of arbitrary Grassmannians \cite{draisma-eggermont}.
\end{remark}

\begin{remark}
The framework behind the theorem above can also be used to give bounds for higher-order syzygies of secant varieties, see \cite{secver-syzygies}. An earlier result for Segre varieties can be found in \cite{deltamod}.
\end{remark}

\section{Twisted commutative algebras}

Both results mentioned above are formal consequences of some algebraic structure being finitely generated. A possible framework for both structures is given by twisted commutative algebras. See \cite{expos} for an introduction.

\subsection{The definition}

Let $A = \bigoplus_{n \ge 0} A_n$ be a $\bZ_{\ge 0}$-graded, associative, unital $A_0$-algebra such that each $A_n$ has a linear action of the symmetric group $S_n$. Embed $S_n \times S_m$ as a subgroup of $S_{n+m}$ by having $S_n$ act on $1,\dots,n$ and $S_m$ act on $n+1,\dots,n+m$ in the natural way. Then $A$ is a \defi{twisted commutative algebra (tca)} if, for all $n$ and $m$,
\begin{itemize}
\item the multiplication map
\[
A_n \otimes A_m \to A_{n+m}
\]
is $S_n \times S_m$-equivariant with respect to the embedding just specified, and 
\item $\tau(xy) = yx$ where $x \in A_n$, $y \in A_m$, and $\tau \in S_{n + m}$ is the permutation that swaps $\{1,\dots,n\}$ and $\{n+1,\dots,n+m\}$ in order, i.e., $\tau(i) = m+i$ if $1 \le i \le n$ and $\tau(n+j) = j$ for $1 \le j \le m$.
\end{itemize}

A basic example is when $E$ is a free $\bk$-module, and $A_n = E^{\otimes n}$ (by convention, $A_0 = \bk$). The action of $S_n$ is by permuting tensor factors, and multiplication is concatenation of tensors. This is the free tca generated in degree $1$ by $E$, and we will denote it by $\Sym(E \langle 1 \rangle)$.

We can also define modules over tca's. Let $A$ be a tca. An $A$-module $M$ is a graded abelian group $M = \bigoplus_{n \ge 0} M_n$ such that each $M_n$ is a linear $S_n$-representation, $M$ is a graded $A$-module in the usual sense, and such that the multiplication map
\[
A_n \otimes M_m \to M_{n+m}
\]
is $S_n \times S_m$-equivariant. A module is \defi{finitely generated} if it can be generated by finitely many elements $m_1, \dots, m_n$ under the action of $A$ and the symmetric groups. We are interested in the category of $A$-modules, and specifically, in formal properties of finitely generated $A$-modules.

\subsection{Categorical model}

For the tca $A = \Sym(E\langle 1 \rangle)$, we can give an alternate model for the category of $A$-modules. For simplicity, suppose that $\dim_\bk E = d< \infty$ and choose a basis $E \cong \bk^d$. 

Define a category $\FI_d$ whose objects are finite sets $S$ and whose morphisms $S \to T$ consist of pairs $(f,g)$ where $f \colon S \to T$ is an injective function and $g \colon T \setminus f(S) \to \{1,\dots,d\}$ is an arbitrary function (``$d$-coloring''). Given another morphism $(f',g') \colon T \to U$, the composition $(f'',g'') \colon S \to U$ is obtained by taking $f'' = f' \circ f$ and $g''(u) = g'(u)$ for $u \in U \setminus f''(S)$. When $d=1$, the function $g$ is superfluous, and $\FI_1 = \FI$ is simply the category of finite sets and injective functions. Note that $\FI_d$ is equivalent to the full subcategory of finite sets of the form $[n] = \{1,\dots,n\}$ (with $[0] = \emptyset$).

Given a commutative ring $\bk$, an \defi{$\FI_d$-module} is a functor $M \colon \FI_d \to \Mod_\bk$. Concretely, we have a $\bk$-module $M(S)$ for every finite set $S$, and $\bk$-linear maps $M(f,g) \colon M(S) \to M(T)$ for every $\FI_d$-morphism $(f,g) \colon S \to T$. Given the remark above, the data of $M$ is completely encoded by the sequence of $\bk$-modules $M_n = M([n])$ together with transition maps $M_n \to M_m$ given by the $\FI_d$-morphisms between $[n]$ and $[m]$. When $S=T$, a morphism is a bijection of the set $S$ to itself, so $M_n$ is a representation of the symmetric group $S_n$.

Given an $\FI_d$-module $M$, define an $A$-module structure on $\bigoplus_{n \ge 0} M_n$ as follows. A basis element of $A_n = E^{\otimes n} \cong (\bk^d)^{\otimes n}$ can be identified with a function $[n] \to [d]$, and we define multiplication by this basis element to be the action of the morphism $(f,g)$ where $f(i) = i+n$ and $g$ is the function $[n] \to [d]$. We omit the (routine) check that this is well-defined.

\begin{proposition}
The functor just defined gives an equivalence between the category of $\FI_d$-modules $($over $\bk)$ and the category of $A$-modules.
\end{proposition}

While the categorical model $\FI_d$ is much easier to work with, we have given the general definition of tca since it suggests a more general algebraic framework. In particular, the category $\FI$ gives an algebraic structure on the cohomology of configuration spaces. 

Let $X$ be a topological space. Given a finite set $S$, let $\Conf_S(X)$ be the set of injective functions $S \to X$, topologized as a subspace of the set of all functions $S \to X$ (a direct product of $|S|$ copies of $X$). Given an injective function $f \colon S \to T$, there is an induced map $\Conf_T(X) \to \Conf_S(X)$ given by composition $(T \to X) \mapsto (S \xrightarrow{f} T \to X)$. By contravariance, for any $i \ge 0$, and coefficient ring $\bk$, one also gets a $\bk$-linear map 
\[
f_* \colon \rH^i(\Conf_S(X); \bk) \to \rH^i(\Conf_T(X); \bk).
\]

An alternative way to phrase this is that each choice of $X, i, \bk$ as above gives a functor $\FI \to \Mod_\bk$ which sends $S$ to $\rH^i(\Conf_S(X); \bk)$ and $f \colon S \to T$ to $f_*$. In the next section, we will see how to exploit this algebraic structure.

\begin{remark}
The category $\FI$ was introduced and studied in \cite{fimodules} over a field of characteristic $0$, and also in a completely different language, in \cite{symc1}. See those papers for more applications.
\end{remark}

\begin{remark}
There is no particular reason to favor symmetric groups over other sequences of groups, such as general linear groups. In fact, ``linear versions'' of tca's are studied in \cite{putman-sam}. These are not spelled out as such, but rather the categorical models (categories ${\bf VI}$, ${\bf VIC}$, and variants) of the simplest cases are the objects of study.
\end{remark}

\subsection{Secant varieties}

The structure behind the example of secant varieties is more involved, so we define the relevant category in the special case where $B = \bk[x_1,\dots,x_r]$ is a polynomial ring. See \cite{secver-syzygies} for more details and motivation.

Define
\[
\multi{r}{d} = \{(x_1, \dots, x_r) \in \bZ_{\ge 0}^r \mid x_1 + \cdots + x_r = d\}.
\]
We think of this as a basis for $\Sym^d \bk^r$, and we define addition componentwise 
\[
+ \colon \multi{r}{d} \times \multi{r}{e} \to \multi{r}{d+e}.
\]

\begin{definition} \label{defn:veronese-cat}
Define the \defi{Veronese category} $\cV_r$ as follows. The objects of $\cV_r$ are pairs $(d,m) \in \bZ_{\ge 0}^2$ and a morphism $\alpha \colon (d,m) \to (e,n)$ consists of the following data:
\begin{itemize}
\item An order-preserving injection $\alpha_1 \colon [m] \to [n]$,
\item A function $\alpha_2 \colon [n] \setminus \alpha_1([m]) \to \multi{r}{e}$,
\item A function $\alpha_3 \colon [m] \to \multi{r}{e-d}$.
\end{itemize}
In particular, $\hom_{\cV_r}((d,m), (e,n)) = \emptyset$ if $d > e$. Given another morphism $\beta \colon (e,n) \to (f,p)$, the composition $\beta \circ \alpha = \gamma \colon (d,m) \to (f,p)$ is defined by 
\begin{itemize}
\item $\gamma_1 = \beta_1 \circ \alpha_1$,
\item $\gamma_2 \colon [p] \setminus \gamma_1([m]) \to \multi{r}{f}$ is defined by:
  \begin{itemize}
  \item if $i \in [p] \setminus \beta_1([n])$, then $\gamma_2(i) = \beta_2(i)$, and
  \item if $i \in \beta_1([n] \setminus \alpha_1([m]))$, then $\gamma_2(i) = \alpha_2(i') + \beta_3(i')$ where $i'$ is the unique preimage of $i$ under $\beta_1$.
  \end{itemize}
\item $\gamma_3 \colon [m] \to \multi{r}{f-d}$ is defined by $\gamma_3(i) = \alpha_3(i) + \beta_3(\alpha_1(i))$. \qedhere
\end{itemize}
\end{definition}

Given $c \in \multi{r}{d}$, let $x^c$ be the corresponding monomial in $x_1,\dots,x_r$. The basic example of a functor on $\cV_r$ is defined by $(d,m) \mapsto (\Sym^d \bk^r)^{\otimes m}$. To a morphism $\alpha \colon (d,m) \to (e,n)$, we define 
\begin{align*}
(\Sym^d \bk^r)^{\otimes m} &\to (\Sym^e \bk^r)^{\otimes n}\\
f_1 \otimes \cdots \otimes f_m &\mapsto g_1 \otimes \cdots \otimes g_n
\end{align*}
where $g_i = f_j x^{\alpha_2(i)}$ if $i = \alpha(j)$, and $g_i = x^{\alpha_2(i)}$ if $i \notin \alpha([m])$.

The keypoint is that this functor encodes certain operations between all of the (non-commutative analogue of the) coordinate rings of the various Veronese embeddings of projective space. One can also work with the actual coordinate rings by introducing a symmetrized version of the Veronese category, but we will avoid discussing this technical point. Fortunately, the ideals of secant varieties (of a fixed secant degree) form a subrepresentation, and the point is to show that this subrepresentation is finitely generated. This will imply the desired result: suppose the subrepresentation is generated in bidegrees $(d_1,m_1), \dots, (d_k,m_k)$. Then for any $e$, the set of elements in bidegrees $(e,m_1), \dots, (e,m_k)$ that come from the generators will be generators for the ideal of the secant variety in the $e$th Veronese embedding; in particular, $\max(m_i)$ is the desired bound in Theorem~\ref{thm:secver}.

\begin{remark}
Modules of the Veronese category do not model modules over a twisted commutative algebra, but there is a relation. The full subcategory on objects of the form $(d,m)$ for $d$ fixed and $m$ varying is the category $\FI_N$ where $N = |\multi{r}{d}|$. Hence, the Veronese category can be interpreted as a certain way of combining many of these categories together.
\end{remark}

\section{Noetherian properties}

We can formalize the situations above as follows. Given a category $\cC$, a $\cC$-module (over a ring $\bk$) is a functor from $\cC$ to the category of (left) $\bk$-modules. In our cases, the isomorphism classes of objects of $\cC$ are indexed by either $\bZ_{\ge 0}$ or $\bZ_{\ge 0}^2$ and so one gets a graded or bigraded sequence of $\bk$-modules. In general, these sequences probably have no interesting properties (or at least no easily usable ones) at all. The key to getting interesting properties is to prove a finite generation statement. 

We say that a $\cC$-module $M$ is \defi{finitely generated} if there exist finitely many elements $m_1, \dots, m_r$ with $m_i \in M(x_i)$ for objects $x_i$ such that each element of $M(y)$ for any object $y$ can be written as a linear combination $\sum_{i=1}^r \alpha_i M(f)(m_i)$ for morphisms $f \colon x_i \to y$. Alternatively, the smallest subfunctor of $M$ that contains the elements $m_1, \dots, m_r$ is $M$. The properties mentioned in \S\ref{sec:conf-eg} are, in fact, properties for general finitely generated $\FI$-modules.

\begin{theorem}
\begin{enumerate}[\rm (1)]
\item Let $M$ be a finitely generated $\FI$-module over a field $\bk$. The function
\[
n \mapsto \dim_\bk M_n
\]
is given by a polynomial function in $n$ for $n \gg 0$.

\item If $\bk$ is a field of characteristic $0$, then $n \mapsto \dim_\bk M_n^{S_n}$ is constant for $n \gg 0$.

\item If $\bk$ is a field of characteristic $p>0$, for each $i$, the function $n \mapsto \dim_\bk \rH^i(S_n; M_n)$ is periodic with period a power of $p$.
\end{enumerate}
\end{theorem}

Property (1) immediately implies Theorem~\ref{thm:CEFN} once we know that the $\FI$-module $S \mapsto \rH^i(\Conf_S(X); \bk)$ is finitely generated. Further given this, property (2) implies the first point of Theorem~\ref{thm:conf-quotient} follows since $\rH^i(\Conf_n(X)/S_n; \bk) = \rH^i(\Conf_n(X); \bk)^{S_n}$. The second point does not quite follow from property (3), instead we refer the reader to \cite{nagpal} for details.

So we see that being finitely generated is a desirable property that we want of the $\FI$-modules that we encounter in \S\ref{sec:conf-eg}. One way to compute these $\FI$-modules is using a version of the Leray spectral sequence (see \cite[\S 4]{CEFN} for details). Under the assumptions made, the initial page of this spectral sequence consists of $\FI$-modules which can be shown directly to be finitely generated. Hence, our desired $\FI$-modules are subquotients of finitely generated ones. Finite generation is clearly inherited by quotient modules, but in complete generality, this need not be inherited by submodules.

We say that a module is \defi{noetherian} if all of its submodules are finitely generated. The following result is exactly what is needed to finish the argument sketched above. It was first proven in the language of tca's over a field of characteristic $0$ in \cite{deltamod}, re-proven (for $d=1$ over characteristic $0$) in \cite{fimodules}, proven in general (for $d=1$) in \cite{CEFN}, and finally in general in \cite{catgb}.

\begin{theorem}
If $\bk$ is a noetherian ring, then every finitely generated $\FI_d$-module over $\bk$ is noetherian.
\end{theorem}

Given a category $\cC$, we will say that the category of $\cC$-modules is \defi{locally noetherian} (over $\bk$) if every finitely generated $\cC$-module over $\bk$ is noetherian. Which $\cC$ are locally noetherian for every noetherian coefficient ring? To finish the discussion from before, we mention the Veronese category.

\begin{theorem}[Sam]
The category of modules of Veronese category over a noetherian ring $\bk$ is locally noetherian. The same is true for the symmetrized version if we restrict to noetherian rings $\bk$ which contain the rational numbers.
\end{theorem}

\subsection{Gr\"obner methods}

A general setup for proving local noetherianity is developed in \cite{catgb} using ideas from Gr\"obner bases and many examples of interest in representation stability can be found there.

We start by recalling a combinatorial proof of the Hilbert basis theorem (i.e., $\bk[x_1, \ldots, x_n]$ is a noetherian ring when $\bk$ is a field) using Gr\"obner bases. Pick an admissible order on the monomials, i.e., a well-order compatible with multiplication. Using the order, we can define initial ideals, and reduce the study of the ascending chain condition to monomial ideals. Now, the set of monomial ideals is naturally in bijection with the set of ideals in the poset $\bN^r$. Thus noetherianity of $\bk[x_1, \ldots, x_n]$ follows from noetherianity of the poset $\bN^r$ (Dickson's lemma), which can be rephrased as the following combinatorial exercise: given infinitely many vectors $v_1, v_2, \ldots$ in $\bN^r$, there exists $i < j$ with $v_i \le v_j$, where $\le$ means coordinate-by-coordinate comparison.

The same idea can be made to work to prove noetherianity of finitely generated $\cC$-modules for various choices of $\cC$. To illustrate this, we will discuss the case of $\FI_d$-modules. There are several reductions that one can make.
\begin{itemize}
\item Define the principal projective $\FI_d$-module generated in degree $n$ to be the functor $P_n$ given by $P_n(S) = \bk[\hom_{\FI_d}([n], S)]$, where $\bk[X]$ means the free $\bk$-module on the set $X$. The action of morphisms is by post-composition. By Yoneda's lemma, a morphism $P_n \to M$ is the same as choosing an element in $M([n])$. Hence, a module $M$ is finitely generated if and only if it is a quotient of $P_{n_1} \oplus \cdots \oplus P_{n_r}$ for some finite $r$.
\item If $M$ and $N$ are noetherian, then the same is true for $M \oplus N$.
\item If $M$ is a quotient of $N$ and $N$ is noetherian, then the same is true for $M$.
\item Putting it all together, to prove that all finitely generated modules are noetherian, it suffices to prove this for $P_n$ for all $n \ge 0$.
\end{itemize}

The $P_n$ are the analogues of free modules, though they behave differently for different $n$ (in the case of graded rings, $P_n$ would simply be a shift of $P_0$, which is why it is sufficient in that case to prove noetherianity of just $P_0$).

To deal with this problem, one would like to reduce to working with {\it monomial subrepresentations} of $P_n$, i.e., those generated by the basis vectors corresponding to morphisms $[n] \to S$. An issue is that there are too few of these for a reduction to be possible: the presence of automorphisms in $\FI_d$ causes many monomials to generate each other. We won't detail this issue any further, but rather discuss how to resolve it. 

Define an auxiliary category $\OI_d$ (O for ordered) as the category of totally ordered sets $S$ where morphisms are pairs $(f,g) \colon S \to T$ as in the definition of $\FI_d$, except now $f \colon S \to T$ is required to be order-preserving. There is a forgetful functor $\Phi \colon \OI_d \to \FI_d$ and we can pullback $\FI_d$-modules to $\OI_d$-modules along $\Phi$. Let $P'_n$ be the corresponding principal projective for $\OI_d$: $P'_n(S) = \bk[\hom_{\OI_d}([n],S)]$. A key point is that $\Phi^*(P_n) \cong (P_n')^{\oplus n!}$, where in the decomposition, an $\FI_d$-morphism $f$ belongs to the copy indexed by the unique permutation $\sigma$ of $[n]$ needed to make $f \circ \sigma$ an order-preserving map. Hence it will actually suffice to show that $P'_n$ is a noetherian $\OI_d$-module.

To do this, one {\it can} reduce to the case of monomial subrepresentations. For this aspect, see \cite[\S 4]{catgb}. To deal with monomial subrepresentations, we can finally reduce to the following combinatorial problem. Put a partial ordering on the set $X_n = \amalg_{m \ge n} \hom_{\OI_d}([n], [m])$ by $\phi \le \phi'$ if there exists $h$ such that $\psi \circ \phi = \phi'$. We want to show that $X_n$ is \defi{well-ordered}, which means that given an infinite sequence $\phi_1, \phi_2, \dots$, there exists $i<j$ such that $\phi_i \le \phi_j$.

We encode morphisms as words in the alphabet $\Sigma = \{*,1,2,\dots,d\}$ as follows. Given a morphism $\phi = (f,g) \colon [n] \to [m]$, create a word $w(\phi)$ of length $m$ whose $i$th letter is a $*$ if $i$ is in the image of $f$, and is $g(i)$ otherwise. Then $\phi$ can be recovered from $w(\phi)$, and furthermore, $\phi \le \phi'$ if and only if $w(\phi)$ is a subsequence of $w(\phi')$.

Given any set $X$, let $X^\star$ be the set of finite words in $X$, with the partial ordering given by the subsequence relation. Higman proved that this is well-ordered whenever $X$ is finite (this can be generalized substantially, but we won't do it here):

\begin{theorem}[Higman's lemma]
If $X$ is finite, then $X^\star$ is well-ordered.
\end{theorem}

By the definition, a subposet of a well-ordered poset is itself also well-ordered, so this solves our problem.

This suggests using a language-theoretic approach to understanding morphisms in $\cC$, and this is explored in \cite[\S 5]{catgb}. We will restrict ourselves to discussing $\OI_1 = \OI$ in this setting. Morphisms are encoded as words in $\{*,1\}$ with $n$ instances of $*$. The set of words that come from a subrepresentation can be shown to be a regular language (and hence the generating function that records the length of words is a rational function), but even further, what we call an {\it ordered language}, which implies that the function counting the number of words of a given length is eventually a polynomial in the length.

\section{What next?}

The big open problem is whether or not finitely generated tca's are noetherian. Some cases of tca's generated in degree $2$ (and their skew-commutative analogues) are proved to be noetherian in \cite{sym2noeth} and \cite{deg2tca} over a field of characteristic $0$ (some related results can be found in \cite{eggermont}). An application of these algebras can be found in \cite{miller-wilson}. 

A categorical model of one case can be given as follows. Let $\cC$ be the category whose objects are sets and morphisms $(f,\Pi) \colon S \to T$ consists of an injection $f \colon S \to T$ and a decomposition $\Pi$ of $T \setminus f(S)$ into $2$-element subsets (so $|S|\equiv |T| \pmod 2$). Then $\cC$-modules are the same as modules over the free tca generated by a trivial element (with respect to the $S_2$-action) of degree 2.

One reason this case is so interesting is that it resists attempts to be solved using Gr\"obner methods. If we follow the outline above, we first define an ordered version $\cC'$ of the category $\cC$ by requiring that the sets are totally ordered and that $f$ is order-preserving. It turns out the module category for $\cC'$ is not locally noetherian! However, pullbacks of $\cC$-modules to $\cC'$ have much more structure which possibly could be utilized. One point here is that almost all known cases of successful Gr\"obner methods applications boil down to Higman's lemma or a variant of it, which might be considered a ``linear'' phenomenon or a ``1-dimensional'' phenomenon. The category $\cC$ is ``non-linear'' in that it deals with graphs rather than words, but it's not clear how to make this more precise.

For any positive integer $d$, modules over the free tca generated by a trivial element of degree $d$ can be modeled as functors on the following category which generalizes $\cC$ above: the objects are finite sets $S$, and a morphism $S \to T$ consists of an injection $f \colon S \to T$ together with a decomposition of $T \setminus f(S)$ as a disjoint union of $d$-element subsets (so $|S| \equiv |T| \pmod d$ if a morphism exists). Almost nothing is known about these module categories when $d \ge 3$. However, over a field of characteristic $0$, this can be transformed into a problem about the $\GL_\infty(\bC)$ action on $\Sym(\Sym^3 \bC^\infty)$, in which case partial progress has been made in \cite{DES}. In particular, \cite{DES} shows that $\Sym(\Sym^3 \bC^\infty)$ is ``topologically noetherian'', i.e., that all $\GL_\infty(\bC)$-equivariant ideals are finitely generated up to radical. The topological noetherianity result was recently extended to $\Sym (F(\bC^\infty))$ for all finite length polynomial functors $F$ by Draisma \cite{draisma}. For some applications of topological noetherianity to uniform bound problems in commutative algebra, see \cite{ESS}.

\end{document}